\documentclass{article}
\usepackage{amsfonts}
\title{\large \bf On class preserving automorphisms of groups of order 32 }
\author{\small \bf Deepak Gumber and Hemant Kalra \\
\small \em School of Mathematics and Computer Applications\\
\small \em Thapar University, Patiala - 147 004,
India\\
\small E-mail: dkgumber@yahoo.com, happykalra26@gmail.com}

\date{}
\newtheorem{thm}{Theorem}[section]
\newtheorem{lm}[thm]{Lemma}
\newtheorem{pp}[thm]{Proposition}

\begin{document}

\maketitle
\begin{abstract}
We study class preserving automorphisms of groups of order thirty two and prove that only two groups have non-inner class preserving automorphisms.
\end{abstract}

\vspace{2ex}


\vspace{2ex}


\section{\large Introduction}

Let $G$ be a finite group and let $p$ be a prime. For $x\in G$, let $x^G$ denote the conjugacy class of $x$ in $G$ and
let $\mathrm{Aut}(G)$ denote the group of all automorphisms of $G$.
An automorphism $\alpha$ of $G$ is called a class preserving
automorphism if $\alpha(x)\in x^G$ for all $x\in G$. Observe that
every inner automorphism of $G$ is a class preserving automorphism.
The set $\mathrm{Aut}_c(G)$ of all class preserving automorphisms is
a normal subgroup of $\mathrm{Aut}(G)$, and the set $\mathrm{Inn}(G)$ of all inner automorphisms is a normal subgroup of
$\mathrm{Aut}_c(G)$. Let $\mathrm{Out}_c(G)$ denote the group
$\mathrm{Aut}_c(G)/\mathrm{Inn}(G)$. In 1911, Burnside [2, p. 463]
posed the question: Does there exist a finite group $G$ which has a
non-inner class preserving automorphism? In other words, whether
there exists a finite group $G$ for which $\mathrm{Out}_c(G)\neq 1$?
 In 1913, Burnside [3] himself gave an
affirmative answer to this question. He constructed a group $G$ of
order $p^6$, $p$ odd, such that $\mathrm{Out}_c(G)\neq 1$. After
this more groups were constructed such that $\mathrm{Out}_c(G)\neq
1$ (see [5], [7], [11], [16] and [17]) and many authors found the
groups for which $\mathrm{Out}_c(G)= 1$ (see [4], [10], [12], [13],
[14]). It follows from [8] that for extra special $p$-groups all
class preserving automorphisms are inner and from [9] it follows
that $\mathrm{Out}_c(G)= 1$ for all groups $G$ of order $p^4$.
Recently Yadav [19] studied the groups of order $p^5$, $p$ an odd
prime. He proved that if $G$ and $H$ are two finite non-abelian
isoclinic groups, then Aut$_{c}(G)\cong$ Aut$_{c}(H)$. He also
proved that Out$_{c}(G)= 1$ for all the groups $G$ except two
isoclinic families. In the present paper we show that out of $51$
groups of order 32, there are only two
groups for which $\mathrm{Out}_c(G)\neq 1$. A list of groups of order 32 is available from Sag and Wamsley's work [15]. We use Sag and Wamsley's list and adopt the same notations for the nomenclature and presentations of the groups. However, we write the generators 1,2,3, and 4 used in [15] respectively as $x,y,z,$ and $w$. All the results proved in section 3 are summarised in the following theorem.\\

\noindent {\bf Theorem A} {\em For all groups G of order $32$, except for the groups
$G_{44}$ and
$G_{45}$, $\mathrm{Out}_c(G)= 1$.}\\

 By $\mathrm{Hom}(G,A)$ we denote the group of all homomorphisms of $G$ into an abelian group $A$. A non-abelian group $G$ that has no non-trivial abelian direct factor is said to be purely non-abelian. For two subgroups $H$ and $K$ of $G$, $[H,K]$ denotes the subgroup of $G$ generated by all commutators $[x,y]=x^{-1}y^{-1}xy$ with $x\in H$ and $y\in K$. By $[x,G]$ we denote the set of all commutators of the form $[x,g],\;g\in G$. Observe that if $G$ is nilpotent of class 2, then $[x,G]$ is a subgroup of $G$. The lower central series of a group $G$ is the descending series $G=\gamma_1(G)\geq\gamma_2(G)\geq\cdots\geq\gamma_i(G)\geq\cdots$, where $\gamma_{n+1}(G)=[\gamma_n(G),G]$, and upper central series is the ascending series $Z(G)=Z_1(G)\leq Z_2(G)\leq\cdots \leq Z_i(G)\leq\cdots$, where $Z_{i+1}(G)=\{x\in G|[x,y]\in Z_i(G)\; \mathrm{for\; all}\; y\in G\}.$ An automorphism $\alpha$ of $G$ is called a central automorphism if $x^{-1}\alpha(x)\in Z(G)$ for each $x\in G$. The set $\mbox{Cent}(G)$ of all central automorphisms of $G$ is a normal subgroup of $\mbox{Aut}(G)$.

\section{\large Toolbox }
All the tools required to prove Theorem A are listed in this section. Fortunately (or unfortunately!) we have been able to prove the theorem using the existing tools. Only one existing tool (Lemma 2.1) is slightly modified (Lemma 2.2).

\begin{lm}
{\em [19, Theorem 3.5]} Let $G$ be a finite $p$-group of class $2$
and let $\{x_1,x_2,\ldots,x_m\}$ be a minimal generating set for $G$
such that $[x_i,G]$  is cyclic, $1\leq i\leq m$. Then
$\mathrm{Out}_c(G)=1$.
\end{lm}
\begin{lm}
If $G$ is a group of order $p^n$, $n\geq 3$ and
$|Z(G)|=p^{n-2}$, then $\mathrm{Out}_c(G)= 1$.
\end{lm}
{\bf Proof.} Observe that nilpotency
class of $G$ is $2$. If $x$ is a non central element
of $G$, then $|C_G(x)|=p^{n-1}$ and thus $|x^G|=|[x,G]|=p$. If $x$ is a central element of $G$,
then $|x^G|=|[x,G]|=1$. In any case, $[x,G]$ is
cyclic and therefore $\mathrm{Out}_c(G)= 1$ by Lemma 2.1.
\hfill $\Box$

\begin{lm}
{\em [9, Proposition 4.1]} Let $G$ be a group which can be generated
by two elements $x, y$ such that every element of $G$ can be written
in the form $x^iy^j$, $i, j$ are integers. Then $\mathrm{Out}_c(G)= 1$.
\end{lm}

\begin{lm}
{\em [8, Proposition 2.2]} Let $G=H\oplus K$ be the direct sum of
its subgroups $H$ and $K$. Then $\mathrm{Out}_c(G)= 1$ if and
only if  $\mathrm{Out}_c(H)=\mathrm{Out}_c(K)=1$.
\end{lm}

\begin{lm}
{\em [6, Proposition 2.7]} Let $G$ be a finite group having an
abelian normal subgroup $A$ with cyclic quotient $G/A$. Then
class-preserving automorphisms of G are inner automorphisms
\end{lm}

\begin{lm}
{\em [19, Lemma 2.2]} Let $G$ be a finite $p$-group such that
$Z(G)\subseteq [x,G]$ for all $x\in G-\gamma_2(G)$. Then
$|\mathrm{Aut}_c(G)|\geq |\mathrm{Cent}(G)||G/Z_{2}(G)|$.
\end{lm}

\begin{lm}
{\em [1]} If $G$ is a purely non abelian finite group, then
$|\mathrm{Cent}(G)|=|\mathrm{Hom}(G/\gamma_{2}(G),Z(G))|$.
\end{lm}

The well known commutator identities
$$[x,yz]=[x,z][x,y][x,y,z];\;\;
[xy,z]=[x,z][x,z,y][y,z],$$
where $x,y,z\in G$, will be frequently used without any reference.

\section{\large Proof of Theorem A}
The groups $G_1$ to $G_7$ are abelian, so we start with non abelian groups. The proof of the theorem follows from all the results proved in this section.

\begin{pp}
For groups $G_{44}$ and $G_{45}$, there exist class preserving automorphisms which are not inner.
\end{pp}
{\bf Proof.} Consider $G=G_{44}$. Since $yx=xy^5$ and $zy=y^7z$, it
follows that $|y|=8$ and every element of $G$ can be written in the
form $x^iy^jz^k$, where $0\leq i,k\leq 1$, and $0\leq j\leq 7$. We
show that $Z(G)=\{1,y^4\}$. Since $[z,y^4]=[z,y^2]^2=[z,y]^4=y^8=1$
and $[x,y^2]=[x,y]^2=y^8=1$,  $y^4\in Z(G)$. It is easy to see that
other powers of $y$ are not in $Z(G)$. Since $xy^4$ does not commute
with $y$, it cannot be in $Z(G)$. If $xy^i\in Z(G);\;i=1,2,3,5,6,7,$
then commuting it with $z$ gives $i\equiv 7i(\mbox{mod}\;8)$, a
contradiction. Using similar arguments, other possibilities can also
be ruled out and thus $Z(G)=\{1,y^4\}$. Consider $G/Z(G)$. Since
$xZ(G)$ commutes with $yZ(G)$ and $zZ(G)$, and $y^2Z(G)$ commutes
with $xZ(G)$ and $zZ(G)$, it follows that $|Z(G/Z(G))|\geq 4$. Since
$yZ(G)\notin Z(G/Z(G)),\;|C_{G/Z(G)}yZ(G)|= 8$ and thus
$|Z(G/Z(G))|=4$. This gives $|Z_2(G)|=8$. Now $G/Z_2(G)$ is abelian,
therefore $G=Z_3(G)$ and hence $G$ is nilpotent of class 3. It is
easy to see that $Z(G)\leq \gamma_2(G)=\Phi(G)$. It follows from
[18, theorem 4.7] that $|\gamma_2(G)|=4$ and  $Z(G)\subseteq [x,G]$
for all $x\in G-\gamma_2(G)$. By Lemma 2.6
$$|\mbox{Aut}_c(G)|\geq |\mbox{Cent}(G)|2^5/|Z_2(G)|.$$
Since $G$ is purely non abelian and $G/\gamma_2(G)$ is elementary abelian of order 8, by Lemma 2.7 we have
$$|\mbox{Cent}(G)|=|\mbox{Hom}(G/\gamma_2(G),Z(G)|=8.$$ Thus  $|\mbox{Aut}_c(G)|\geq 2^5>2^4=|\mbox{Inn}(G)|$. That $\mathrm{Out}_c(G_{45})\neq 1$ can be shown similarly. \hfill $\Box$\\

\noindent {\bf Remark:} The order of $\mbox{Aut}_c(G)$ in the above
proposition is in fact exactly equal to 32. Observe that
$C_G(x)=\;<x,y^2,z>,\;C_G(y)=\;<y>$ and $C_G(z)=\;<x,y^4,z>$. Thus
$|C_G(x)|=16$ and $|C_G(y)|=|C_G(z)|=8$. Hence $|x^G|=2$ and
$|y^G|=|z^G|=4$. Since any class preserving automorphism preserves
the conjugacy classes, there are $|x^G|,\;|y^G|$ and $|z^G|$ choices
for the images of  $x,\;y$ and $z$ respectively under it. Thus
$|\mbox{Aut}_c(G)|\leq |x^G||y^G||z^G|=32$.

\begin{pp}
The groups $G_{42}$ and $G_{43}$ are extraspecial and therefore all class preserving automorphisms are inner.
\end{pp}
{\bf Proof.} Consider $G=G_{42}$. Since $1=[y^2,w]=[x^2,w]=[x,w]^2=x^4$, it follows that $|x|=|y|=|z|=|w|=4$ and $x^2\in Z(G)$. Since $wx=xw^3$ and $zy=y^3z$, every element of $G$ can be written as $x^iy^jz^kw^l$, where $0\leq i\leq 3$ and $0\leq j,k,l\leq 1$. If
$z(x^iy)= (x^iy)z$, then $x^iy^3z=x^iyz$ and thus $y^2=1$ which is not so. So $x^iy\notin Z(G)$. Using similar arguments we can show that $x^iz,x^iw,x^iyz,x^iyw,x^izw$, and $x^iyzw$ cannot be in the center. Thus $Z(G)=\{x^2,1\}$. Observe that any element of $G/Z(G)$ is of the form $x^iy^jz^kw^lZ(G)$, where $0\leq i,j,k,l\leq 1$. Since $xw=wx^3\equiv wx\;(\mbox{mod}\;Z(G))$ and $yz=zy^3\equiv zy\;(\mbox{mod}\;Z(G))$, it follows that $G/Z(G)$ is abelian and hence $\gamma_2(G)=Z(G)=\Phi(G)$. Similar arguments show that $G_{43}$ is also an extraspecial $p$-group.
 \hfill $\Box$

\begin{pp}
If $G$ is any of the groups from $G_{23}$ to $G_{25}$, then $G=H\oplus K$ for some subgroups $H$ and $K$  such that $\mathrm{Out}_c(H)=\mathrm{Out}_c(K)= 1$ and hence $\mathrm{Out}_c(G)= 1$ by Lemma $2.4$.
\end{pp}
{\bf Proof.} Observe that $G_{23}=H\oplus K$, where $H=\{x,y|
x^8=y^2=1, yx=x^7y\}$ and $K=\{z|z^2=1\}$; $G_{24}=H\oplus K$, where
$H=\{x,y| x^8=y^2=1, yx=x^3y\}$ and $K=\{z| z^2=1\}$; and
$G_{25}=H\oplus K$, where $H=\{x,y| x^8=1, y^2=x^4,
y^{-1}xy=x^{-1}$\} and $K=\{z| z^2=1\}$. Since $K$ is cyclic  and every element of $H$ can be written as $x^iy^j$ for suitable $i$ and $j$, $\mathrm{Out}_c(K)=\mathrm{Out}_c(H)= 1$.
 \hfill $\Box$

\begin{pp}
If $G$ is any of the groups from $G_{22}$, $G_{29}$, $G_{30}$,
$G_{32}$, and $G_{49}-G_{51}$, then every element of $G$ can be
written as $x^iy^j$, for suitable integers $i$ and$j$, and hence
$\mathrm{Out}_c(G)= 1$ by Lemma $2.3$.
\end{pp}
{\bf Proof:} Observe that all the groups are generated by $x$ and $y$, and $yx$ can be written as $xy^9,x^7y,x^3y,x^7y,x^{15}y,x^7y$ and $x^{15}y$ in the respective groups.
 \hfill $\Box$

\begin{pp}
If $G$ is any of the groups from $G_{8}$ to $G_{21}$, then
$|Z(G)|=8$ and hence $\mathrm{Out}_c(G)= 1$ by Lemma $2.2$.
\end{pp}
{\bf Proof.} For $G_8$ it is clear that $z$ and $w$ are in the
center. Since $[x^2,y]=[x,y]^2=1, x^2$ is in the center and thus
$|Z(G_8)|=8$. For $G_9$ again it is obvious that $z$ and $w$ are in
the center. Since $x^4=[y,x]^2=[y^2,x]=1$, $|x|=4$ and since
$[x^2,y]=[x,y]^2=1, x^2$ is in the center and thus $|Z(G_9)|=8$. For
$G_{10}$, $x$ and $w$ are in the center and since
$x^4=[z,y]^2=[z^2,y]=1, |x|=4$ and therefore $|Z(G_{10})|=8$. For
$G_{11}$, clearly $z$ is in the center and since
$[x,y^2]=[x,y]^2=[x^2,y]=1$, it follows that $y^2$ is in the center
and $|[x,y]|=2$. Also, since $z$ commutes with $x$ and $y$, it
commutes with $[x,y]$ as well.  Thus $[x,y]$ is also in the center
and hence $|Z(G_{11})|=8$. For $G_{12}$, $z$ is in the center and
since $[x^2,y]=[x,y^2]=1$, both $x^2$ and $y^2$ are in the center
and thus $|Z(G_{12})|=8$. Observe that every element of $G_{13}$ can
be written in the form $x^iy^jz^k$; $0\leq i\leq 1$, $0\leq j\leq
7$, and $0\leq k\leq 1$. Therefore $|y|=8$ and $|x|=2$. Now clearly
$z$ is in the center and $[y^2,x]=[y,x]^2=1$ implies that $y^2$ is
also in the center. Therefore $|Z(G_{13})|=8$. For $G_{14}$, clearly
$z$ is in the center and since $[x^2,y]=[x,y]^2=1$, $x^2$ is in the
center and thus $|Z(G_{14})|=8$. In $G_{15}$, observe that
$x^4=[y,x]^2=[y^2,x]=1=[y,x^2]$. Thus $|x|=4$ and $x^2$ is in the
center. Also $z$ is in the center and therefore $|Z(G_{15})|=8$. In
$G_{16}$, $|x|=4$ and $x$ and $z^2$ are in the center. Thus
$|Z(G_{16})|=8$. In $G_{17}$, $x$ is in the center and $|x|=8$. Thus
$|Z(G_{17})|=8$. In $G_{18}$, $[x,y]$ is in the center. Also
$[x^2,y]=[x,y^2]=[x,y]^2=1$ implies that $x^2$ and $y^2$ are in the
center and $|[x,y]|=2$. Therefore $|Z(G_{18})|=8$. In $G_{19}$,
since $[x^2,y]=[x,y^2]=[x,y]^2=1$, it follows that $x^2$ and $y^2$
are in the center and $|Z(G_{19})|=8$. In $G_{20}$ clearly $[x,y]$
is in the center and $[x,y^2]=[x,y]^2=[x^2,y]=1$ implies that
$|[x,y]|=2$ and $y^2$ is in the center. Thus $|Z(G_{20})|=8$. In
$G_{21}$ it follows from the relation  $[x^2,y]=[x,y^2]=1$ that
$x^2$ and $y^2$ are in the center and $|Z(G_{21})|=8$. \hfill $\Box$

\begin{pp}
If $G$ is any of the groups from $G_{26}-G_{28}$, $G_{31}$, $G_{33}-G_{41}$ and $G_{46}-G_{48}$, then $G$ has an abelian normal subgroup $H$ such that $G/H$ is cyclic and hence $\mathrm{Out}_c(G)= 1$ by Lemma $2.5$.
\end{pp}
{\bf Proof.} In $G_{26}$, $1=[x^2,z]=[y^4,z]=[y^2,z]^2=[y,z]^4$.
Thus $|x|=4,|y|=8$, and $H=\;<x,y>.$ In $G_{27}$,
$1=[x,y^2]=[x,y]^4=x^8$ implies that $|x|=8$. Let $z=[x,y]$, then
$[z,x]=1$ and $|z|=4$. Thus $H=\;<x,z>$. In $G_{28}$,
$1=[x^4,y]=[x^2,y]^2=[x,y]^4=x^8$ implies that $|x|=8$. Take
$z=[x,y]$, then $[z,x]=1,|z|=4$ and $H=\;<x,z>$. In $G_{31}$, take
$z=[x,y]$, then $z$ commutes with $x$. Also,
$1=[x^4,y]=[x^2,y]^2=[x,y]^4=z^4$ implies that $|z|=4$ and thus
$H=\;<x,z>$. In $G_{33}$, observe that each of the elements
$x,y,zxz^{-1},$ and $zyz^{-1}$ is of order 2 and commutes with
another. Thus $H=\;<x,y,zxz^{-1},zyz^{-1}>$.
 In $G_{34}$, it is easily seen that $H=\;<x,y>$. In $G_{35}$, observe that $|x|=4$ and $H=\;<x,y>$.
In $G_{36}$ and $G_{37}$, $1=[x,z]^2=[x^2,z]$ and $[x,z]=x^2$ implies that $xz=zx^3$. Take $w=[y,z]$, then $|w|=2$ and it commutes with $x$ and $y$. Thus $H=\;<x,y,w>$. In $G_{38}$, $1=[x^2,z]=[x,z]^2=y^4\;;\;(y^2z)^2=([z,x]z)^2=1$. Thus $|y|=4$ and $y^2z=zy^2$. Take $w=[y,z]$, then $|w|=2$ and it commutes with $y$ and $z$. Thus $H=\;<x,y,w>$. In $G_{39}$, $[x^2,z]=[x,z]^2=x^4\; ;\;1=[x,z^2]=[x,z]^2[x^2,z]=x^8\; ;\;[y^2,z]=[y,z]^2=x^4\; ;\;1=[y,z^2]=[y,z]^2[y,z,z]=x^4[x^2y^2,z]=x^4[x^2,z][y^2,z]
=x^{12}.$ Thus $|x|=4$ and $H=\;<x,y>$.
In $G_{40}$, $z^4=[z,y]^2=[z^2,y]=[x^{-2}y^2,y]=1\; ;\;[y^2,z]=[y,z]^2=1\; ;\;x^4=[x,z]^2=[x^2,z]=[y^2z^2,z]=1.$ Thus $|x|=4,|y|=4$ and $H=<x,y>$.
 In $G_{41}$, $[x^2,z]=[x,z]^2=x^4y^4\; ;\;[y^2,z]=[y,z]^2=x^4\; ;\;1=[x,z^2]=[x,z]^2[x^2y^2,z]=x^{12}y^{8}\; ;\;
 1=[y,z^2]=[y,z]^2[x^2,z]=x^8y^4.$ Thus $|x|=4,|y|=4$ and
$H=<x,y>$. In $G_{46}$, let $z=[x,y]$ and $w=[x,y,y]=[z,y]$, then
$z$ commutes with $x$ and $w$ commutes with $y$. Observe that
$1=[x^2,y]=[x,y]^2=z^2\;;\;[x,y^2]=[x,y]^2[x,y,y]=w\;
;\;1=[x,y^4]=[x,y^2]^2[x,y^2,y^2]=w^2.$ Thus $|z|=|w|=2$.  Since
$xy^2=(xy)y=(yxz)y=yx(zy)=yx(yzw)=y(xy)zw=y(yxz)zw=y^2xw$,
$wx=(xy^2xy^2)x=(xy^2x)y^2x=(xy^2x)(xwy^2)=xw.$ Thus $w$ commutes
with $x$ and hence with $z$ as well. It is easy to see that every
element of $G_{46}$ can be written in the form $x^iy^jz^kw^l$, where
$0\leq i,k,l\leq 1$ and $0\leq j\leq 3$. Let $H=\;<x,z,w>$, then $H$
is abelian, $Hy=yH$ and $G_{46}/H=\;<yH>.$ In $G_{47}$ let
$[x,y]=z$, then $z$ commutes with $x$ and $[z,y]=y^4.$ Observe that
$1=[x^2,y]=[x,y]^2=z^2\;;\;[z,y^2]=[z,y]^2=y^8\;;\;$ and
$1=[z^2,y]=[z,y][z,y,z][z,y]=y^4[y^4,z]y^4=y^4[y^2,z]^2y^4= y^{-8}.$
Thus $y^2$ commutes with $z$, $|z|=2$ and $|y|=8$. Since
 $zy=y^5z$ and $yx=xyz$, every
 element of $G_{47}$ can be written as $x^iy^jz^k$, where
 $0\leq i,k\leq 1$ and $0\leq j\leq 7$. Let $H=\;<x,y^4,z>,$ then it is easy to verify that $H$ is a normal subgroup of $G_{47}$ and $G_{47}/H=<Hy>$
 is cyclic. In
 $G_{48}$, let $[x,y]=z$, then $z$ commutes with $x$ and $z^2=[x^2,y]=1$. Since $1= [z^2,y]=[z,y]^2=y^8$, $|x|=4$ and $|y|=8$.
 Now $zy=y^5z$ and $yx=xyz$, therefore every element of $G_{48}$ can be written as
 $x^iy^jz^k$, where $0\leq i,k\leq 1$ and $0\leq j\leq 7$. Let
 $H=\;<x,z>$, then $Hy=yH$ and $G_{48}/H=\;<Hy>$.
\hfill $\Box$

\vspace{.2in}
 \noindent {\bf Acknowledgement:} The research of the
second author is supported by Council of Scientific and Industrial
Research, Government of India. The same is gratefully acknowledged.


\begin{thebibliography}{20}
\bibitem{}J. E. Adney and T. Yen, {\em Automorphisms of a $p$-group}
Illinois J. Math. {\bf 9} (1965), 137-143.
\bibitem{}W. Burnside, {\em Theory of groups of finite order}, 2nd
Ed. Dover Publication, Inc., 1955. Reprint of the 2nd edition
(Cambridge, 1911).
\bibitem{}W. Burnside, {\em On the outer automorphisms of a group},
Proc. London Math. Soc. (2) {\bf 11} (1913), 40-42.
\bibitem{}M. Fuma and Y. Ninomiya, {\em``Hasse principle" for finite
p-groups with cyclic subgroups of  index $p^2$}, Math. J. Okayama
Univ. {\bf 46} (2004), 31-38.
\bibitem{}H. Heineken, {\em Nilpotente Gruppen, deren s$\ddot{a}$mtliche Normalteiler charakteristisch sind}, Arch.
Math. (Basel) {\bf 33} (1980), No. 6, 497-503.
\bibitem{}M. Hertweck and E. Jespers, {\em Class-preserving automorphisms and normalizer property for
Blackburn groups}, J. Group Theory {\bf 12} (2009), 157-169.
\bibitem{}D. Jonah and M. Konvisser, {\em Some non-abelian p-groups with abelian automorphism
groups}, Arch. Math. (Basel) {\bf 26} (1975), 131-133.
\bibitem{}M. Kumar and L. R. Vermani,
{\em ``Hasse principle" for extraspecial p-groups}, Proc. Japan
Acad. {\bf 76}, Ser. A, No.8 (2000), 123-125.
\bibitem{}M. Kumar and L. R. Vermani, {\em ``Hasse principle" for the groups
of order $p^4$}, Proc. Japan Acad. {\bf 77}, Ser. A, No. 6 (2001),
95-98.
\bibitem{}M. Kumar and L. R. Vermani, {\em On automorphisms of some
p-groups}, Proc. Japan Acad. {\bf 78}, Ser. A, No. 4 (2002), 46-50.
\bibitem{}I. Malinowska, {\em On quasi-inner automorphisms of a finite p-group}, Publ. Math. Debrecen {\bf
41}
(1992), No. 1-2, 73-77.
\bibitem{}T. Ono, {\em ``Hasse principle" for $GL_2(D)$}, Proc.
Japan Acad. {\bf 75} Ser. A (1999), 141-142.
\bibitem{}T. Ono and H. Wada, {\em ``Hasse principle" for free groups}, Proc.
Japan Acad. {\bf 75} Ser. A (1999), 1-2.
\bibitem{}T. Ono and H. Wada, {\em ``Hasse principle for symmetric and
alternating groups}, Proc. Japan Acad. {\bf 75} Ser. A (1999),
61-62.
\bibitem{}T. W. Sag and J. W. wamsley, {\em Minimal presentations
for groups of order $2^n$, $n\leq 6$}, J. Aust. Math. Soc. {\bf 15}
(1973), 461-469.
\bibitem{}C. H. Sah, {\em Automorphisms of finite groups}, Journal of Algebra {\bf 10} (1968), 47-68.
\bibitem{}F. Szechtman, {\em n-inner automorphisms of finite groups}, Proc. Amer. Math. Soc. {\bf 131} (2003),
3657-3664.
\bibitem{}H. Tendra and W. Moran, {\em Flatness conditions on finite
$p$-groups}, Comm. Algebra {\bf 32} (2004), 2215-2224.
\bibitem{}M. K. Yadav, {\em On automorphisms of some finite p-groups}, Proc.
Indian Acad. Sci. (Math Sci.) {\bf 118} (1) (2008), 1-11.


\end{thebibliography}
\end{document}